\documentclass[a4paper, leqno, 11pt]{article}
\usepackage{latexsym}
\usepackage{amsmath}
\usepackage{amssymb}
\usepackage{enumerate}
\usepackage{theorem}
\usepackage{array}
\newtheorem{theorem}{Theorem}[section]
\newtheorem{lemma}[theorem]{Lemma}
\newtheorem{corollary}[theorem]{Corollary}
\newtheorem{proposition}[theorem]{Proposition}
\newtheorem{atheorem}{Theorem A.\!\!}
\newtheorem{alemma}[atheorem]{Lemma A.\!\!}

\theorembodyfont{\rmfamily}

\newtheorem{setup}[theorem]{Setup}
\newtheorem{remark}[theorem]{Remark}

\newcommand{\proof}{\noindent \mbox{\em Proof.\hspace*{2mm}}}
\newcommand{\qed}{\hfill \mbox{$  \Box $}}
\newcommand{\ssgyokan}{\vskip 8pt}
\newcommand{\sgyokan}{\vskip 10pt}
\newcommand{\gyokan}{\vskip 20pt}
\title{Nonsingular plane filling curves of minimum degree over a finite field
and their automorphism groups:
Supplements to a work of Tallini}
\author{Masaaki Homma
\thanks{Partially supported by Grant-in-Aid
for Scientific Research (19540058), JSPS.}
\\
 Department of Mathematics,
Kanagawa University\\
Yokohama 221-8686, Japan\\
homma@n.kanagawa-u.ac.jp
\and Seon Jeong Kim
\thanks{Partially supported by 
the Korea Research Foundation Grant funded by the
Korean Government(MOEHRD) (KRF-2006-312-C00016).}\\
Department of Mathematics and RINS\\
Gyeongsang National University\\
Jinju 660-701, Korea \\
skim@gnu.kr }
\date{}

\begin{document}

\maketitle
\begin{abstract}
Our concern is a nonsingular plane curve defined over
a finite field ${\Bbb F}_q$
which includes all the ${\Bbb F}_q$-rational points
of the projective plane.
The possible degree of such a curve is at least $q+2$.
We prove that
nonsingular plane curves of degree $q+2$
having the property actually exist.
More precisely, we write down explicitly
all of those curves.
Actually, Giuseppe Tallini studied such curves in his old paper in 1961.
We explain the connection between his work and ours.
Moreover we give another proof of his result on the automorphism group
of such a curve, from the viewpoint of linear algebra.
\end{abstract}

\section{Introduction}
To start with, we clarify the setting of our concern.
\begin{setup}
Let $q$ be a power of a prime number,
and ${\Bbb F}_q$ the finite field consisting of $q$ elements.
The algebraic closure of ${\Bbb F}_q$ is denoted by $K$.
We consider the projective plane ${\Bbb P}^2$ defined over ${\Bbb F}_q$
and denote by ${\Bbb P}^2({\Bbb F}_q)$
the set of ${\Bbb F}_q$-rational points of ${\Bbb P}^2$.
A plane curve over ${\Bbb F}_q$ is said to be nonsingular
if the curve is nonsingular at any $K$-point, not only
the ${\Bbb F}_q$-points.
\end{setup}
Under this setup, 
an interesting question is whether
there exists a nonsingular plane curve over ${\Bbb F}_q$
including ${\Bbb P}^2({\Bbb F}_q)$,
which is the simplest case of a question posed by
N.~Katz \cite[Question~10]{kat}.
The question of Katz was settled affirmatively
by O.~Gabber \cite{gab} and B.~Poonen \cite{poo} independently.
So we may ask
what the smallest degree of such a curve is.

If $C$ is a nonsingular plane curve defined over ${\Bbb F}_q$
such that $ C \supset {\Bbb P}^2({\Bbb F}_q)$,
then the tangent line $l$ at a point $P \in {\Bbb P}^2({\Bbb F}_q)$
to $C$ meets with $C$ at least $q$ points other than $P$
because $l$ is also defined over ${\Bbb F}_q$.
Moreover, since the intersection multiplicity of $l$ and $C$
at $P$ is at least $2$,
the degree of $C$ is at least $q+2$ by B\'{e}zout's theorem.
Let $ x, y, z $ be a system of homogeneous coordinates of
${\Bbb P}^2$ over ${\Bbb F}_q$.
We denote by ${\frak h}$
the homogeneous ideal of the set ${\Bbb P}^2({\Bbb F}_q)$
in ${\Bbb F}_q[x,y,z]$,
and by ${\frak h}_d$ its homogeneous part of degree $d$.
Then any element of ${\frak h}_{q+2}$ can be represented as
\[F_A = (x,y,z) A \, {}^{t}(y^q z - yz^q, z^q x-  z x^q, x^q y -x y^q),\]
where $A \in GL(3, {\Bbb F}_q)$.
In Section 3, we show that {\em the curve $F_A=0$ is nonsingular
if and only if the characteristic polynomial of $A$ is irreducible
over ${\Bbb F}_q$}.

Half a century ago,
G. Tallini
considered a problem which seems to be close to ours
\cite{tal1, tal2}
\footnote{
A summary of these works can be found in
\cite[Chap. 8, Exercise~12]{hir-kor-tor}.
}.
Namely he considered an irreducible curve instead of a nonsingular curve
in our context.
We explain relations between his results and ours in Section~4.

After we wrote up the earlier version on this topic,
an anonymous reviewer pointed out that
the irreducibility of a plane filling curve of degree $q+2$,
which was proved by Tallini,
would imply the smoothness,
if one used the latter part of \cite{tal2}
on the automorphism group of that curve.
Certainly, it is true, however Tallini's study of automorphisms of such a curve
involves the study of certain invariants $\theta_{ijk}$ of the curve.
In Section 5, we give another proof of his results with a correction
in the context of linear algebra,
and explain the remark by the anonymous reviewer.
In the appendix, we study the difference between
the image of 
the center of $M \in GL(n, {\Bbb F}_q)$ in $PGL(n, {\Bbb F}_q)$
and the center of the image of $M$ in $PGL(n, {\Bbb F}_q)$
when the characteristic polynomial of $M$ is irreducible,
which is necessary for Section~5.
\section{The homogeneous ideal of plane filling curves}
Although natural generators of  ${\frak h}$ (Prop.~\ref{generatorofh})
and the property of ${\frak h}_{q+1}$ (Prop.~\ref{ideal})
are already known by Tallini \cite{tal1, tal2},
we give their proofs for reader's convenience.
\begin{proposition}[Tallini]\label{generatorofh}
The homogeneous ideal ${\frak h}$ is generated by
\[
(y^q z - yz^q, z^q x-  z x^q, x^q y -x y^q)
\]
over ${\Bbb F}_q[x,y,z]$.
\end{proposition}
\proof
Denote by
$
{\frak i} = ( y^q z - yz^q, z^q x - z x^q, x^q y -x y^q)
$
Then it is obvious that
$
{\frak i} \subseteq {\frak h}.
$
Let $f(x,y,z) \in {\frak h}_d$.
Write it as
\[
f(x,y,z) = z g(x,y,z) + h(x,y),
\]
where $g(x,y,z)$ is homogeneous of degree $d-1$
and $h(x,y)$ is a homogeneous polynomial in $x, y$
of degree $d$.
Since $h(x,y) = f(x,y,0)$,
it vanishes on ${\Bbb P}^1({\Bbb F}_q)$ with coordinates $x$ and $y$.
Hence
$xy \prod_{\alpha \in {\Bbb F}_q^{\times}}(x-\alpha y)$
divides $h(x,y)$.
So it is enough to show that
$zg(x,y,z) \in {\frak i}$
because
$
xy \prod_{\alpha \in {\Bbb F}_q^{\times}}(x-\alpha y)
 = x^q y - x y^q.
$
To see this, it is sufficient to see that
$ g(x,y,1) \in (y^q -y , x-x^q)$
in ${\Bbb F}_q[x,y]$.
Note that $ g(\alpha ,\beta ,1) =0$
for all $(\alpha , \beta ) \in ({\Bbb F}_q)^2$.
Denote by $I = (y^q -y , x-x^q)$.
Since $ x^q \equiv x \mod I$ and $ y^q \equiv y \mod I$,
there is a $q \times q$ matrix $M$ with entries in ${\Bbb F}_q$
such that
\[
g(x,y,1) \equiv
(x^{q-1}, x^{q-2}, \ldots , 1) M
\left(
\begin{array}{c}
y^{q-1}\\
y^{q-2}\\
\vdots \\
1
\end{array}
\right)
\mod I.
\]
Hence
\[
\left(
\begin{array}{cccc}
&\vdots& & \\
\alpha^{q-1}&\alpha^{q-2}&\cdots&1\\
&\vdots& & 
\end{array}
\right)_{\alpha \in {\Bbb F}_q}
M
\left(
\begin{array}{ccc}
&\beta^{q-1}&\\
\cdots & \beta^{q-2}& \cdots \\
       & \vdots & \\
       & 1&
\end{array}
\right)_{\beta \in {\Bbb F}_q}
= 0.
\]
Since the first and the third matrices in the above
are invertible,
we have $M=0$. Hence $g(x,y,1) \in I$.
\qed

\ssgyokan

We denote by
\[
\left\{
  \begin{array}{ccl}
  U &=& y^q z - y z^q \\
  V &=& z^q x - z x^q \\
  W &=& x^q y - x y^q .
  \end{array}
\right.
\]
We observe the behavior of $U, V$ and $W$
under a linear transformation
\begin{equation}\label{lineartransform}
 \left(
  \begin{array}{c}
  x\\
  y\\
  z
  \end{array}
 \right)
  = B \left(
  \begin{array}{c}
  x'\\
  y'\\
  z'
  \end{array}
 \right)
\end{equation}
with
\[
B = 
\left(
  \begin{array}{ccc}
  b_{11} & b_{12} & b_{13} \\
  b_{21} & b_{22} & b_{23} \\
  b_{31} & b_{32} & b_{33} 
  \end{array}
\right)
\in
GL(3, {\Bbb F}_q).
\]
We need more notation:
\[
\left\{
  \begin{array}{ccl}
  U' &=& {y'}^q {z'} - {y'} {z'}^q \\
  V' &=& {z'}^q {x'} - {z'} {x'}^q \\
  W' &=& {x'}^q {y'} - {x'} {y'}^q 
  \end{array}
\right.
\]
and $\tilde{b}_{ij}$ denotes the $(i,j)$-cofactor of $B$,
for example,
$\tilde{b}_{13} = (-1)^{1+3}
    \det \left(
       \begin{array}{cc}
        b_{21} & b_{22} \\
        b_{31} & b_{32}
       \end{array}
    \right).$
Hence
\[
\left(
  \begin{array}{ccc}
  \tilde{b}_{11} & \tilde{b}_{12} & \tilde{b}_{13} \\
  \tilde{b}_{21} & \tilde{b}_{22} & \tilde{b}_{23} \\
  \tilde{b}_{31} & \tilde{b}_{32} & \tilde{b}_{33} 
  \end{array}
\right)
= (\det B){}^{t}\!B^{-1}
\]
and $B^{(q)} = B$, where $B^{(q)}=( b_{ij}^q )$.
The following lemma can be proved easily,
but is essential to our linear-algebraic point of view.
\begin{lemma}\label{UVW}
Under the above notation, we have
\[
 \left(
  \begin{array}{c}
  U \\
  V \\
  W
  \end{array}
 \right)
  = (\det B) {}^{t}\!B^{-1} 
 \left(
  \begin{array}{c}
  U' \\
  V' \\
  W'
  \end{array}
 \right).
\]
\end{lemma}
\proof
By straightforward computations, we have
\begin{eqnarray*}
U &=& (\begin{array}{ccc} x& y & z \end{array})
      \left(
       \begin{array}{ccc}
       0 & 0 &0 \\
       0 & 0 & -1 \\
       0 & 1 & 0
       \end{array}
      \right)
         \left(
         \begin{array}{c} 
         x^q \\
         y^q \\
         z^q 
         \end{array}
         \right)  \\
    &=& (\begin{array}{ccc} x'& y' & z' \end{array})
        {}^{t}B  
      \left(
       \begin{array}{ccc}
       0 & 0 &0 \\
       0 & 0 & -1 \\
       0 & 1 & 0
       \end{array}
      \right)
       B^{(q)}
       \left(
         \begin{array}{c} 
         {x'}^q \\
         {y'}^q \\
         {z'}^q 
         \end{array}
         \right) \\
     &=& (\begin{array}{ccc} x'& y' & z' \end{array})
        \left(
       \begin{array}{ccc}
       0 & -\tilde{b}_{13} & \tilde{b}_{12} \\
       \tilde{b}_{13} & 0 & -\tilde{b}_{11} \\
        -\tilde{b}_{12}& \tilde{b}_{11}& 0
       \end{array}
      \right) 
        \left(
         \begin{array}{c} 
         {x'}^q \\
         {y'}^q \\
         {z'}^q 
         \end{array}
         \right)
      \\
     &=& \tilde{b}_{11} U' + \tilde{b}_{12} V' + \tilde{b}_{13}W' ;\\
   V &=& \tilde{b}_{21} U' + \tilde{b}_{22} V' + \tilde{b}_{23}W' ;\\
   W &=& \tilde{b}_{31} U' + \tilde{b}_{32} V' + \tilde{b}_{33}W' ,
\end{eqnarray*}
which is the desired formula.
\qed

\ssgyokan

\begin{proposition}[Tallini]\label{ideal}
Let $F=a_1 U+ a_2 V + a_3W$ be a nonzero element of ${\frak h}_{q+1}$,
where $a_1, a_2, a_3 \in {\Bbb F}_q$.
Then the curve defined by $F=0$ has a unique singular point
$(a_1, a_2, a_3)$,
and is the union of $q+1$ ${\Bbb F}_q$-lines passing through the point.
\end{proposition}
\proof
Since
\begin{eqnarray*}
F_x &=& a_2 z^q - a_3 y^q = (a_2 z - a_3 y)^q \\
F_y &=& a_3 x^q - a_1 z^q = (a_3 x - a_1 z)^q \\
F_z &=& a_1 y^q - a_2 x^q = (a_1 y - a_2 x)^q ,
\end{eqnarray*}
only the solution of
$F_x = F_y = F_z = 0$
in ${\Bbb P}^2$ is $(a_1, a_2, a_3)$.
Choose a matrix $B \in GL(3, {\Bbb F}_q)$ so that
$
{}^{t}\!(a_1, a_2, a_3)
 = 
B {}^{t}\!(0,0,1).$
Using new coordinates $x', y' , z'$
with
(\ref{lineartransform}) and Lemma~\ref{UVW},
the curve is given by
$W'=0$.
The curve is obviously the union of $q+1$ ${\Bbb F}_q$-lines
passing through $( 0 , 0 , 1 )$ in the new coordinates.
\qed
\section{The condition of smoothness for a member of ${\frak h}_{q+2}$}
From Prop.~\ref{generatorofh},
any member of ${\frak h}_{q+2}$ can be written as
$
l_1U + l_2 V + l_3 W ,
$
where
$
l_j = a_{1j}x +  a_{2j}y +  a_{3j}z
$
is a linear form over ${\Bbb F}_q$
for each $j = 1,2,3$.
In other words, it takes the form of
\[
F_A :=
\left(
 \begin{array}{ccc}
  x &
  y &
  z
 \end{array}
\right)
A
\left(
  \begin{array}{c}
  U \\
  V \\
  W
 \end{array}
  \right),
\]
where
$A = ( a_{ij} )$ is a $3\times 3$ matrix whose entries are in
${\Bbb F}_q$.
Note that $F_A$ may represent a null form even if $A$ is a nonzero matrix.

\begin{lemma}\label{null}
For an element $F_A$ of ${\frak h}_{q+2}$,
$F_A=0$ as an element of ${\Bbb F}_q[x,y,z]$
if and only if
$A = \mu E$ $(\mu \in {\Bbb F}_q).$
\end{lemma}
\proof
The {\it if} part is obvious.
If $F_A(x,y,z)$ is $0$ as a polynomial,
then so are $F_A(x,y,0)$, $F_A(x,0,z)$ and $F_A(0,y,z)$.
Hence $a_{13} = a_{23} = 0$, $a_{12} = a_{32} = 0$ and
$a_{21} = a_{31} = 0$, and hence
\[
F_A(x,y,z) = a_{11}x (y^qz - yz^q) + a_{22}y(z^qx - zx^q)
               + a_{33}z(x^qy - xy^q).
\]
This polynomial represents $0$ only if
$a_{11} = a_{22} = a_{33}$
\qed

\ssgyokan

We denote by $C_A$ the curve defined by
$F_A =0$
in ${\Bbb P}^2$,
and by $f_A(t)$
the characteristic polynomial
$\det (tE - A)$
of $A$,
where $E$ is the unit matrix of degree $3$.
Since $\deg f_A(t) =3$,
it is irreducible over ${\Bbb F}_q$
if and only if no eigen-value of $A$ is in ${\Bbb F}_q$.
\begin{theorem}\label{mainth}
The curve $C_A$ is nonsingular if and only if
the characteristic polynomial $f_A(t)$ is irreducible over ${\Bbb F}_q$.
\end{theorem}
\proof
We denote by
$
F = l_1U + l_2 V + l_3 W.
$
Then
\begin{equation}\label{partialFgen}
 \left\{
  \begin{array}{rcl}
  F_x &=& a_{11} U + a_{12} V + a_{13} W
        + l_2 z^q - l_3 y^q \\
  F_y &=& a_{21} U + a_{22} V + a_{23} W
        + l_3 x^q - l_1 z^q \\
  F_z &=& a_{31} U + a_{32} V + a_{33} W
        + l_1 y^q - l_2 x^q .
  \end{array}
 \right.
\end{equation}

The first claim is that
{\em $C_A$ is nonsingular at any ${\Bbb F}_q$-point if and only if
no eigen-value of $A$ is in ${\Bbb F}_q$.}
For an ${\Bbb F}_q$-point $(\alpha, \beta, \gamma)$,
it is a solution of $F_x=F_y=F_z=0$ if and only if
a solution of
\begin{equation}\label{solpartialFgen}
l_2 z - l_3 y = l_3 x - l_1 z = l_1 y - l_2 x =0
\end{equation}
because of (\ref{partialFgen}) and of the identities
$\alpha^q = \alpha$, $\beta^q = \beta$ and $\gamma^q = \gamma$.
The last condition (\ref{solpartialFgen}) means that
\begin{equation}\label{solpartialFgen2}
\frac{a_{11}\alpha + a_{21}\beta + a_{31}\gamma}{\alpha}
= 
\frac{a_{12}\alpha + a_{22}\beta + a_{32}\gamma}{\beta}
=
\frac{a_{13}\alpha + a_{23}\beta + a_{33}\gamma}{\gamma}
\end{equation}
holds, that is,
${}^{t}\! A
{}^{t}\!(\alpha , \beta , \gamma )
  = \lambda {}^{t}\!(\alpha , \beta , \gamma )$
for some $\lambda$.
Since all quantities appeared in (\ref{solpartialFgen2}) are in ${\Bbb F}_q$,
so is the quantity $\lambda$, which is an eigen-value of $A$.
Therefore the first claim has been proved.

\ssgyokan

The second claim is that
{\em if no eigen-value of $A$ is in ${\Bbb F}_q$,
then $C_A$ is nonsingular at any $K$-point.}
The proof of this claim is divided into two steps.
In the first step, we reduce the polynomial $F$ to a simpler form.
In the second step, we prove the curve to be nonsingular by using
an idea of St\"{o}hr and Voloch \cite{sto-vol}.

(Step 1)
Since any eigen-value of $A$ is not contained in ${\Bbb F}_q$,
the characteristic polynomial $f_A(t)$,
say $t^3 -(ct^2 + bt + a)$, is irreducible over ${\Bbb F}_q$.
Note that the characteristic polynomial of
\begin{equation}\label{canonicalform}
A_0 = 
\left(
\begin{array}{ccc}
0&0&a\\
1&0&b\\
0&1&c
\end{array}
\right)
\end{equation}
is also $f_A(t)$.
Since $f_A(t)$ is irreducible over ${\Bbb F}_q$,
we know there exists $B \in GL(3, {\Bbb F}_q)$ such that
$
{}^{t}\!B A{}^{t}\!B^{-1} = A_0
$
by a standard linear algebra (e.g. Gantmacher \cite[VI, \S 3]{gan}).
Choose a new system of coordinates $(x', y', z')$ as
$(x, y, z) = (x', y', z'){}^{t}\!B$.
Then
\begin{eqnarray*}
F &=& \left(
 \begin{array}{ccc}
  x &
  y &
  z
 \end{array}
\right)
 A
\left(
  \begin{array}{c}
  U \\
  V \\
  W
 \end{array}
  \right) \\
&=& 
(\det B)
\left(
 \begin{array}{ccc}
  x' &
  y' &
  z'
 \end{array}
\right)
{}^{t}\!B
A
{}^{t}\!B^{-1}
 \left(
  \begin{array}{c}
  U' \\
  V' \\
  W'
  \end{array}
 \right)
\mbox{\rm \ (by Lemma~\ref{UVW})}\\
&=&
(\det B) 
\left(
 \begin{array}{ccc}
  x' &
  y' &
  z'
 \end{array}
\right)
A_0
 \left(
  \begin{array}{c}
  U' \\
  V' \\
  W'
  \end{array}
 \right)\\
&=&
(\det B)
\left(
y'U' +z'V'+(ax'+by'+cz')W'
\right).
\end{eqnarray*}
So we may start from the following setting.
Our curve $C$ is defined by $F=0$, where
\begin{equation}\label{canform}
F=  y(y^q z - y z^q)
  + z(z^q x - z x^q)
  + (ax+by+cz)(x^q y - xy^q)
\end{equation}
and
\begin{equation}\label{charpoly}
t^3 - (ct^2 + bt +a)
\end{equation}
is irreducible over ${\Bbb F}_q$, which is the characteristic polynomial of
$A_0$.
Note that from the first claim,
$C$ is nonsingular at any ${\Bbb F}_q$-point of ${\Bbb P}^2$.

(Step 2)
If the curve $C$ has a singular point $R$,
it is a solution of
$F_x = F_y = F_z = 0$,
where
\begin{equation}\label{partialF}
\left\{
\begin{array}{rcl}
F_x &=& a(x^qy-xy^q) + z^{q+1} -(ax+by+cz)y^q \\
F_y &=& (y^qz -yz^q) + b(x^qy-xy^q) -yz^q  + (ax+by+cz)x^q \\
F_z &=& (z^qx-zx^q) + c(x^qy-xy^q) +y^{q+1} - zx^q
\end{array}
\right.
\end{equation}
Hence it is also a solution of
$
G = x^q F_x + y^q F_y + z^q F_z = 0.
$
Now we consider another curve $D$ defined by $G=0$.
Then $R \in C \cap D$.
The polynomial can be expressed as
\[
G = y^q(y^qz -yz^q) + z^q(z^qx-zx^q)
 + (ax+by+cz)^q (x^qy-xy^q)
\]
and is a member of ${\frak h}_{2q+1}$.
Since 
\begin{equation}\label{partialG}
\left\{
\begin{array}{rcl}
G_x & = & \left( z^2 - (ax+by+cz)y \right)^q \\
G_y & = & \left( (ax+by+cz)x - yz \right)^q \\
G_z & = & (y^2 - zx)^q
\end{array}
\right. ,
\end{equation}
the solutions of $G_x = G_y = G_z = 0$ in ${\Bbb P}^2$ are
\[
\{
Q_{\lambda} = (\lambda^{-2}, \lambda^{-1}, 1)
\mid
\mbox{\rm $\lambda$ is a root of $t^3 - (ct^2 + bt +a) = 0$}
\}.
\]
Since $G(\lambda^{-2}, \lambda^{-1}, 1)=0$,
the set of singular points of $D$ consists of those three points.
Moreover, we have
$F(\lambda^{-2}, \lambda^{-1}, 1)=0$
by direct computation.
Hence
$
C \cap D \supseteq
\{ Q_{\lambda} \} \cup {\Bbb P}^2({\Bbb F}_q).
$

Next we estimate the intersection number
$i(C.D;P)$ for $P \in {\Bbb P}^2({\Bbb F}_q)$
and $i(C.D;Q_{\lambda})$.
For a moment, we suppose $C$ and $D$ to have no common component.
First we consider
$P = (\alpha, \beta, \gamma ) \in {\Bbb P}^2$
with $\alpha, \beta, \gamma \in  {\Bbb F}_q$,
which is nonsingular on both curves $C$ and $D$.
Then by (\ref{partialF}) and (\ref{partialG}),
we have
\[
\begin{array}{ccccl}
F_x(\alpha, \beta, \gamma ) &=& G_x(\alpha, \beta, \gamma )
&=& \gamma^2 -(a \alpha + b \beta + c\gamma)\beta \\
F_y(\alpha, \beta, \gamma ) &=& G_y(\alpha, \beta, \gamma )
&=& (a \alpha + b \beta + c\gamma)\alpha - \beta \gamma \\
F_z(\alpha, \beta, \gamma ) &=& G_z(\alpha, \beta, \gamma )
&=& \beta^2 - \alpha \gamma,
\end{array}
\]
which means that $C$ and $D$ have the common tangent line at $P$.
Hence $i(C.D;P) \geq 2$.
Secondly, we consider behavior of $C$ and $D$
around $Q_{\lambda} =(\lambda^{-2}, \lambda^{-1}, 1)$.
We choose local coordinates around $Q_{\lambda} \in {\Bbb A}^2$
as
\[
\left\{
\begin{array}{ccc}
s & = & x - \lambda^{-2} \\
t & = & y - \lambda^{-1}.
\end{array}
\right.
\]
Hence local equations of $C$ and $D$ around $Q_{\lambda}$
are given by 
$F(s+ \lambda^{-2},  t+ \lambda^{-1}, 1)=0$
and 
$G(s+ \lambda^{-2},  t+ \lambda^{-1}, 1)=0$
respectively.
Write as
\[
F(s+ \lambda^{-2},  t+ \lambda^{-1}, 1)
= (\mbox{\rm Coeff}_s F)s + (\mbox{\rm Coeff}_t F)t
  + \tilde{F}(s,t),
\]
where $\mbox{\rm Coeff}_s F, \mbox{\rm Coeff}_t F \in K$
and $\deg \tilde{F}(s,t) \geq 2$.
To compute $\mbox{\rm Coeff}_s F$ and  $\mbox{\rm Coeff}_t F$,
recall the equation $a \lambda^{-2} + b \lambda^{-1} +c = \lambda$.
Then
\begin{eqnarray}
\mbox{\rm Coeff}_s F &=& F_x(\lambda^{-2}, \lambda^{-1},1) \nonumber \\
&=& (\lambda^{1-q} -1)(a \lambda^{-q-2} -1),
\label{coeffsF}
\end{eqnarray}
and
\begin{eqnarray}
\mbox{\rm Coeff}_t F &=& F_y(\lambda^{-2}, \lambda^{-1},1) \nonumber \\
&=& (\lambda^{1-q} -1)(b \lambda^{-q-2} +\lambda^{-q} + 2\lambda^{-1}).
\label{coefftF}
\end{eqnarray}
Hence
$\mbox{\rm Coeff}_s F \neq 0$.
In fact, since $\lambda \not\in {\Bbb F}_q$,
$\lambda^{q-1} \neq 1$.
Moreover, since $\lambda$ is a root of $t^3 - (ct^2 + bt +a)=0$
which is irreducible over ${\Bbb F}_q$,
other roots of the cubic equation are
$\lambda^{q}$ and $\lambda^{q^2}$.
Hence $a = \lambda^{q^2+q+1}$.
So $a \lambda^{-q-2} -1 = \lambda^{q^2-1} -1 \neq 0$
because $\lambda \not\in {\Bbb F}_{q^2}$ either.
Therefore $C$ is nonsingular at $Q_{\lambda}$ and
the tangent line to $C$ at $Q_{\lambda}$ is given by
\begin{equation}\label{tanCatQ}
(\lambda^{1-q} -1)(a \lambda^{-q-2} -1)s
 + (\lambda^{1-q} -1)(b \lambda^{-q-2} +\lambda^{-q} + 2\lambda^{-1})t
     = 0.
\end{equation}
For $D$, put
\[
G(s+ \lambda^{-2},  t+ \lambda^{-1}, 1)
= (\mbox{\rm Coeff}_{s^q} G)s^q + (\mbox{\rm Coeff}_{t^q} G)t^q
  + \tilde{G}(s,t),
\]
where $\deg \tilde{G}(s,t) \geq q+1$.
Here
\begin{eqnarray}
\mbox{\rm Coeff}_{s^q} G &=& a^q(\lambda^{-2q-1}- \lambda^{-q-2}) 
                          +\lambda^{q-1} - 1  \nonumber \\
&=& (\lambda^{1-q} -1)(a \lambda^{-q-2} -\lambda^{q-1})
  \mbox{\rm  \ \ (because $a^q = a$)},
\label{coeffsG}
\end{eqnarray}
and
\begin{eqnarray}
\mbox{\rm Coeff}_{t^q} G  &=& b^q(\lambda^{-2q-1} - \lambda^{-q-2}) 
        - \lambda^{q-2} + ( \lambda^{-q} - \lambda^{-1} )
               + \lambda^{-q} \nonumber \\
&=& (\lambda^{1-q} -1)(b \lambda^{-q-2} +\lambda^{q-2} + 2\lambda^{-1}).
\label{coefftG}
\end{eqnarray}
By the similar arguments to those in proving $\mbox{\rm Coeff}_{s} F \neq 0$,
we know $\mbox{\rm Coeff}_{s^q} G \neq 0$.
In particular, the multiplicity $\mu_{Q_{\lambda}}(D)$
of $D$ at $Q_{\lambda}$ is $q$.
So $i(C.D;Q_{\lambda}) \geq q$.
Summing up, we have
\begin{eqnarray*}
(\deg C)(\deg D) &=& (C.D) \\
     &=& \sum_{\lambda} i(C.D; Q_{\lambda})
          + \sum_{P \in {\Bbb P}^2({\Bbb F}_q)} i(C.D; P) \\
       &\geq& 3q + 2(q^2+q+1)
\end{eqnarray*}
because $C$ and $D$ have no common component.
Since $(\deg C)(\deg D) =(q+2)(2q+1) = 3q + 2(q^2+q+1)$,
equality holds above,
which implies $C \cap D = \{ Q_{\lambda} \}_{\lambda} \cup {\Bbb P}^2({\Bbb F}_q)$.
So if there is a singular point of $C$,
then it must be one of the $Q_{\lambda}$'s.
However $C$ is nonsingular at $Q_{\lambda}$
as we have already seen.
So we can conclude that $C$ is nonsingular.

The remaining task is to show that $C$ and $D$ have no common component.
If they have a common component $E$, then there is another component $E'$
of $D$ because $\deg D > \deg C$.
Since any point of $E \cap E'$ is a singular point of $D$,
$E$ contains one of the $Q_{\lambda}$'s
because only those three points are singularities of $D$.
Since $C$ is nonsingular at $Q_{\lambda}$,
$E$ must coincide with $C$ around $Q_{\lambda}$.
Therefore the tangent line (\ref{tanCatQ}) to $C$ at $Q_{\lambda}$
must be contained in the tangent cone
\[
\left(
(\mbox{\rm Coeff}_{s^q} G)^{1/q}s + (\mbox{\rm Coeff}_{t^q} G)^{1/q}t
\right)^q = 0
\]
of $D$ at $Q_{\lambda}$.
Therefore if
\[
\det
\left(
\begin{array}{cc}
(\mbox{\rm Coeff}_{s} F)^q & (\mbox{\rm Coeff}_{t} F)^q \\
\mbox{\rm Coeff}_{s^q} G   & \mbox{\rm Coeff}_{t^q} G
\end{array}
\right)
\neq 0,
\]
then we can conclude that
$C$ and $D$ have no common component.

Now we compute the determinant,
in which we use the notation
$|M|$ instead of $\det M$.
It is equal to
\[
\begin{array}{rl}
&
 \left|
  \begin{array}{cc}
  (\lambda^{1-q} -1)^q (a \lambda^{-q-2} -1)^q &
     (\lambda^{1-q} -1)^q(b \lambda^{-q-2} +\lambda^{-q} + 2\lambda^{-1})^q \\
  (\lambda^{1-q} -1)(a \lambda^{-q-2} -\lambda^{q-1}) &
        (\lambda^{1-q} -1)(b \lambda^{-q-2} +\lambda^{q-2} + 2\lambda^{-1})
  \end{array}
\right|
\\
=
&
(\lambda^{1-q} -1)^{q+1}
 \left|
  \begin{array}{cc}
     (a \lambda^{-q-2} -1)^q &
         (b \lambda^{-q-2} +\lambda^{-q} + 2\lambda^{-1})^q \\
     a \lambda^{-q-2} -\lambda^{q-1} &
       b \lambda^{-q-2} +\lambda^{q-2} + 2\lambda^{-1}
  \end{array}
\right|
\end{array}
\]
by (\ref{coeffsF}), (\ref{coefftF}), (\ref{coeffsG}) and (\ref{coefftG}).
By straightforward computation,
we have
\[
\begin{array}{rl}
&  \left|
  \begin{array}{cc}
     (a \lambda^{-q-2} -1)^q &
         (b \lambda^{-q-2} +\lambda^{-q} + 2\lambda^{-1})^q \\
     a \lambda^{-q-2} -\lambda^{q-1} &
       b \lambda^{-q-2} +\lambda^{q-2} + 2\lambda^{-1}
  \end{array}
\right| \\
=&
  a  \left|
  \begin{array}{cc}
     \lambda^{-q^2-2q}  &
         \lambda^{-q^2} + 2\lambda^{-q} \\
     \lambda^{-q-2}  &
       \lambda^{q-2} + 2\lambda^{-1}
  \end{array}
\right|
   +b  \left|
  \begin{array}{cc}
          -1 & \lambda^{-q^2-2q} \\
          -\lambda^{q-1} & \lambda^{-q-2}
  \end{array}
\right|
   +  \left|
  \begin{array}{cc}
      -1 &
         \lambda^{-q^2} + 2\lambda^{-q} \\
      -\lambda^{q-1} &
          \lambda^{q-2} + 2\lambda^{-1}
  \end{array}
\right| \\
=&
 (1- \lambda^{q^2-1})\lambda^{-q^2-1}(2a\lambda^{-2}+b\lambda^{-1}+\lambda)^q.
\end{array}
\] 
If this is $0$, then
$2a\lambda^{-2}+b\lambda^{-1}+\lambda$ must be 0.
Hence
$\lambda^3 +b \lambda + 2a =0$.
However, since
$
\lambda^3 -(c\lambda^2 + b \lambda +a)=0,
$
\[
c\lambda^2 + 2b \lambda + 3a =0.
\]
If $c$ or $2b$ is nonzero, then $\lambda$ is a root of
a polynomial of degree less than $3$, which is a contradiction.
Therefore $c=2b =0$.
If the characteristic is not $3$, this implies $a=0$, which is
a contradiction because $a = \lambda^{1+q +q^2} \neq 0$.
If the characteristic is $3$,  $c=2b =0$ implies that
the minimal polynomial of $\lambda$ is
$t^3 - a$, which is absurd.
This completes the proof.
\qed

\begin{corollary}
The smallest degree of a nonsingular plane curve over
${\Bbb F}_q$ which contains ${\Bbb P}^2({\Bbb F}_q)$
is $q+2$.
\end{corollary}
\proof
Since the degree of such a curve is at least $q+2$,
it is enough to show the existence of a nonsingular member of degree
$q+2$ in $\frak h$.
To do this,
choose an irreducible polynomial (\ref{charpoly})
over ${\Bbb F}_q$
and consider the curve defined by (\ref{canform}).
\qed

\section{Connection with a work of Tallini}
As was mentioned in the Introduction,
Tallini proved the following theorem \cite{tal2}.
Our coordinates $x,y,z$ correspond to $x_1, x_0, x_2$
in \cite{tal1, tal2} in this order.

\begin{theorem}[Tallini]
Any irreducible member of ${\frak h}_{q+2}$
is projectively equivalent over ${\Bbb F}_q$ to 
an element in the form {\rm (\ref{canform})}
such that the polynomial {\rm (\ref{charpoly})}
is irreducible,
and vice versa.
\end{theorem}
Essentially he proved that
{\em
a member $F$ of ${\frak h}_{q+2}$
is irreducible over $K$ if and only if
the curve $C$ defined by $F=0$ is
nonsingular at each ${\Bbb F}_q$-point,
}
which seems a combinatorial nature according to his proof.
On the other hand, in the proof of Th.~\ref{mainth},
we have seen that the similar statement
just replaced the word {\em irreducible} by {\em nonsingular}
holds true.
So we have
\begin{corollary}
For the curve $C_A$ defined by an element $F_A$ of ${\frak h}_{q+2}$,
the following conditions are equivalent{\rm :}
\begin{enumerate}[{\rm (a)}]
\item $C_A$ is nonsingular{\rm ;}
\item $C_A$ is irreducible over $K${\rm ;}
\item $C_A$ is nonsingular at each ${\Bbb F}_q$-point{\rm ;}
\item the characteristic polynomial $f_A(t)$ of $A$ is irreducible.
\end{enumerate}
\end{corollary}

Moreover, Tallini gave a classification of those curves as follows.
\begin{theorem}[Tallini]
Any irreducible member of ${\frak h}_{q+2}$ is projectively
equivalent to one of the following forms over ${\Bbb F}_q${\rm :}
\begin{enumerate}[{\rm (i)}]
\item $yU+zV+a(x+y)W$, where
 $t^3-at-a$ is irreducible over ${\Bbb F}_q${\rm ;}
\item $yU+zV+axW$, where
 $t^3-a$ is irreducible over ${\Bbb F}_q$,
 which happens only in the case $q \equiv 1 \mod 3${\rm ;}
\item $q = 3^e$ and $yU+zV-(x+ az)W$,
where $t^3 + at^2 + 1$ is irreducible over ${\Bbb F}_q$.
\end{enumerate}
\end{theorem}

Tallini's classification can be understood in our context as follows.
\begin{theorem}
For two nonsingular curves $C_A$ and $C_B$,
they are projectively equivalent over ${\Bbb F}_q$ if and only if
there are $\rho \in {\Bbb F}_q^{\times}$
and $\mu \in  {\Bbb F}_q$
such that
\begin{equation}\label{equivrelcubics}
f_A(t) = \rho^3 f_B(\frac{t-\mu}{\rho}).
\end{equation}
\end{theorem}
\proof
This comes from Lemma~\ref{UVW} and Lemma~\ref{null}
with the property of characteristic polynomials:
$f_{\mu E + \rho B}(t) = \rho^3 f_{B}(\frac{t-\mu}{\rho})$.
\qed

\ssgyokan

Taking account of Theorem~\ref{mainth},
we know that Tallini's list of the classification corresponds
to a set of complete representatives
of irreducible cubics in $t$ over ${\Bbb F}_q$
under the equivalence relation (\ref{equivrelcubics}).

\section{Automorphism groups of nonsingular plane filling curves
of degree $q+2$}
The purpose of this section is to study the automorphism group
${\rm Aut}_{{\Bbb F}_q}(C_A)$ of $C_A$ over ${\Bbb F}_q$.
Since the smoothness of $C_A$ is already established in Section~3,
any automorphism comes from a linear transformation of ${\Bbb P}^2$,
because $g^2_{q+2}$ is unique
\cite[Appendix~A, Exercises 17 and 18]{arb-cor-gri-har}.
So we can regard ${\rm Aut}_{{\Bbb F}_q}(C_A)$ as a subgroup of
$PGL(3, {\Bbb F}_q)$.
Let $GL(3, {\Bbb F}_q) \to PGL(3, {\Bbb F}_q)$
be the natural homomorphism,
and a bar over an object in $GL(3, {\Bbb F}_q)$ 
denotes its image by this map.
Let $Z_{GL(3, {\Bbb F}_q)}({}^t\!A)$
be the center of ${}^t\!A \in GL(3, {\Bbb F}_q)$,
and $Z_{PGL(3, {\Bbb F}_q)}(\overline{{}^t\!A})$
the center of $\overline{{}^t\!A} \in PGL(3, {\Bbb F}_q)$.
Then we have
\begin{equation}\label{center}
\overline{Z_{GL(3, {\Bbb F}_q)}({}^t\!A)}
\subset Z_{PGL(3, {\Bbb F}_q)}(\overline{{}^t\!A})
\subset {\rm Aut}_{{\Bbb F}_q}(C_A)
\end{equation}
by Lemma~\ref{UVW}.

Let $f_A(t) = t^3 -(ct^2 + bt +a)$
be the characteristic polynomial of $A$
and $\{ \lambda, \lambda^{q}, \lambda^{q^2} \}$
the roots of $f_A(t) =0$.
Let $\Lambda_0$ be an eigen-vector of ${}^t\!A$ with the eigen-value
$\lambda$,
that is ${}^t\!A\Lambda_0 = \lambda \Lambda_0$.
Then $\Lambda_0^{(q^i)}$, say $\Lambda_i$, is an eigen-vector of 
${}^t\!A$ with the eigen-value
$\lambda^{q^i}$ for $i = 0, 1, 2$.
We denote by $\overline{\Lambda}_i$
the point of ${\Bbb P}^2$ corresponding to
the column vector $\Lambda_i$ for $i = 0, 1, 2$.
These three points agree with the $Q_{\lambda}$'s
in the proof of Theorem~\ref{mainth}
if we choose $A$ as (\ref{canonicalform}).
\begin{lemma}\label{propertiesB}
Let $\overline{B} \in  {\rm Aut}_{{\Bbb F}_q}(C_A)$.
Then,
\begin{enumerate}[{\rm (a)}]
\item there are $\rho = \rho_B \in {\Bbb F}_q^{\times}$
and $\mu = \mu_B \in {\Bbb F}_q$
such that
$
{}^t\! A B = \rho B {}^t\! A + \mu B {\rm; \ } and
$
\item $\{ \overline{B}\overline{\Lambda}_0,
          \overline{B}\overline{\Lambda}_1,
          \overline{B}\overline{\Lambda}_2 \}
        = \{
        \overline{\Lambda}_0, \overline{\Lambda}_1,
        \overline{\Lambda}_2
        \}$.
\end{enumerate}
\end{lemma}
\proof
Since $\overline{B}(C_A)$ is defined by
\[
(x,y,z){}^t\!BA{}^t\!B^{-1}
\left(
  \begin{array}{c}
  U \\
  V \\
  W
 \end{array}
  \right)
 = 0,
\]
$\overline{B}(C_A) = C_A$ if and only if
${}^t\!BA{}^t\!B^{-1} - \rho A = \mu E$
for some $\rho \in {\Bbb F}_q^{\times}$
and $\mu \in {\Bbb F}_q$
by Lemma~\ref{null}.
This completes the proof of (a).

From (a), we have
\begin{eqnarray*}
{}^t\!AB \Lambda_i &=& \rho B{}^t\!A \Lambda_i
                        + \mu B \Lambda_i \\
                   &=& (\rho \lambda^{q^i} + \mu)B \Lambda_i.
\end{eqnarray*}
So $B \Lambda_i$ is an eigen-vector of ${}^t\!A$.
Hence $\overline{B} \overline{\Lambda}_i 
= \overline{\Lambda}_{\sigma_B(i)}$ for some
$\sigma_B(i) \in \{0,1,2\}.$
\qed

\ssgyokan

From this lemma, we can define the group homomorphism
${\rm Aut}_{{\Bbb F}_q}(C_A) \to S_3$ by
$\overline{B} \mapsto \sigma_B$.
We denote by $\pi$ this group homomorphism.

\begin{lemma}\label{kerandim}
\begin{enumerate}[{\rm (a)}]
\item For $\overline{B} \in {\rm Aut}_{{\Bbb F}_q}(C_A)$,
the following conditions are equivalent{\rm :}
 \begin{enumerate}[{\rm (i)}]
   \item $\overline{B} \overline{\Lambda}_i =\overline{\Lambda}_i$
     for some $i \in \{0,1,2 \}${\rm ;}
   \item $\overline{B} \in \ker \pi${\rm ;}
   \item $\overline{B} \in
          \overline{Z_{GL(3, {\Bbb F}_q)}({}^t\!A)}$.
 \end{enumerate}
In particular, $\ker \pi = \overline{Z_{GL(3, {\Bbb F}_q)}({}^t\!A)}.$
\item If $\pi$ is nontrivial, then ${\rm Im}\, \pi = A_3$.
\end{enumerate}
\end{lemma}
\proof
(a)\  Suppose $\overline{B} \overline{\Lambda}_i =\overline{\Lambda}_i$
for some $i$.
Then $B \Lambda_i = \kappa \Lambda_i$
for some $\kappa \in K^{\times}$.
Since
${}^t\!AB \Lambda_i = \rho B{}^t\!A \Lambda_i + \mu B \Lambda_i$
by Lemma~\ref{propertiesB}(a),
$\kappa \lambda^{q^i} \Lambda_i
   = \kappa (\rho \lambda^{q^i} + \mu ) \Lambda_i.$
So $\lambda^{q^i} = \rho \lambda^{q^i} + \mu$
with $\rho,  \mu \in {\Bbb F}_q$.
Hence $\rho = 1$ and $\mu = 0$,
which means $B \in Z_{GL(3, {\Bbb F}_q)}({}^t\!A)$.

Next suppose $B \in Z_{GL(3, {\Bbb F}_q)}({}^t\!A)$.
Then, for $j=0, 1, 2$,
${}^t\!AB \Lambda_j = B {}^t\!A \Lambda_j
 = \lambda^{q^j} B\Lambda_j$,
which means $B\Lambda_j$ is an eigen-vector of ${}^t\!A$
with the eigen-value $\lambda^{q^j}$.
So $\overline{B} \overline{\Lambda}_j = \overline{\Lambda}_j$,
which means $\overline{B} \in \ker \pi$.
The implication ${\rm (ii)} \Rightarrow {\rm (i)}$
is obvious.

(b)\  From  ${\rm (i)} \Rightarrow {\rm (ii)}$ of (a),
${\rm Im}\, \pi$ does not contain any transposition.
\qed

\begin{lemma}\label{replace}
For $C_A$, there is a matrix $A' \in GL(3, {\Bbb F}_q)$
such that $C_{A'} = C_A$ and
${\rm Aut}_{{\Bbb F}_q}(C_A)=Z_{PGL(3, {\Bbb F}_q)}({}^t\!\overline{A'}),$
except the case where $q=3^e$ and the characteristic polynomial
$f_A(t)$ of $A$ is of the form
$f_A(t)= t^3 -(bt + a)$
with a square element $b=\mu^2$ $(\mu \in {\Bbb F}_q)$.
In the exceptional case,
$Z_{PGL(3, {\Bbb F}_q)}({}^t\!\overline{A})
= \overline{Z_{GL(3, {\Bbb F}_q)}({}^t\!A)}$ and
${\rm Aut}_{{\Bbb F}_q}(C_A)/Z_{PGL(3, {\Bbb F}_q)}({}^t\!\overline{A})
\simeq A_3$.
\end{lemma}
\proof
From Lemma~\ref{kerandim} (a) with (\ref{center}),
the assertion is obvious if ${\rm Im}\, \pi$ is trivial.
So we suppose ${\rm Im}\, \pi =A_3$.
Choose $\overline{B}_1\in {\rm Aut}_{{\Bbb F}_q}(C_A)$
so that $\pi (\overline{B}_1)$ is not the identity.
By Lemma~\ref{propertiesB} (a),
there are $\rho_1, \mu_1 \in {\Bbb F}_q^{\times}$
such that
$
{}^t\!A B_1 = \rho_1 B_1 {}^t\!A + \mu_1 B_1.
$
Taking the trace of the both sides of
\begin{equation}\label{eq_b1}
B_1^{-1} {}^t\!A B_1 = \rho_1 {}^t\!A + \mu_1 E,
\end{equation}
we have
$(1- \rho_1) {\rm tr}\, A = 3 \mu_1.$
When $q \neq 3^e$, $\rho_1=1$ implies
$\overline{B}_1 \in \ker \pi$ by Lemma~\ref{kerandim},
which contradicts with the assumption on $\overline{B}_1$.
Hence $\rho_1 \neq 1$.
Put $A' = A + \frac{\mu_1}{\rho_1 - 1}E.$
Then $C_{A'} = C_A$ by Lemma~\ref{null}, and
\begin{eqnarray*}
B_1^{-1} {}^t\!A' B_1 &=& B_1^{-1} {}^t\!A B_1 +\frac{\mu_1}{\rho_1 - 1}E\\
      &=& \rho_1 {}^t\!A + (\mu_1 + \frac{\mu_1}{\rho_1 - 1})E
        = \rho_1 {}^t\!A' ,
\end{eqnarray*}
which means $\overline{B}_1 \in Z_{PGL(3, {\Bbb F}_q)}({}^t\!\overline{A'}).$
For a given $\overline{B} \in {\rm Aut}_{{\Bbb F}_q}(C_{A'})$,
choose an integer $s = 0$ or $1$ or $2$
so that $\pi(\overline{B}\overline{B_1^{-s}})$ is the identity.
Hence $\overline{B}\overline{B^{-s}_1} \in 
    \overline{Z_{GL(3, {\Bbb F}_q)}({}^t\!A')}.$
Hence $\overline{B} \in 
  \overline{B}_1^s \cdot \overline{Z_{GL(3, {\Bbb F}_q)}({}^t\!A')}
   \subset Z_{PGL(3, {\Bbb F}_q)}({}^t\!\overline{A'}).$
   
When $q= 3^e$, we should do more carefully.
If $\rho \neq 1$ in (\ref{eq_b1}),
the same arguments work well even if $q= 3^e$.
So we have to consider the case where
\[
B_1^{-1} {}^t\!A B_1 =  {}^t\!A + \mu_1 E 
           \   \   \  \mbox{\rm  with $\mu_1 \in {\Bbb F}_q$}
\]
holds.
Comparing the sets of eigen-values of both sides above,
we have
\[
\{  \lambda ,  \lambda^{q}, \lambda^{q^2} \}
 = \{  \lambda + \mu ,  \lambda^{q} + \mu, \lambda^{q^2}+ \mu \}.
\]
If $\lambda + \mu = \lambda$, then
$B_1 \in Z_{GL(3, {\Bbb F}_q)}({}^t\!A')$,
which contradicts to the choice of $B_1$.
Hence $\lambda + \mu = \lambda^q$ or $\lambda^{q^2}$.
In either case,
we have
\[
\{  \lambda ,  \lambda^{q}, \lambda^{q^2} \}
 = \{ \lambda, \lambda + \mu , \lambda + 2\mu  (= \lambda - \mu) \}.
\]
Taking the norm of $\lambda$ over ${\Bbb F}_q$, we have
$ \lambda(\lambda +\mu)(\lambda - \mu) =a \in {\Bbb F}_q^{\times}.$
So $\lambda$ satisfies the equation
$t^3 - (\mu^2 t +a) =0,$
which must be the characteristic polynomial of $A$.

For the remaining statement, it is enough to see
that there is a matrix $B \in GL(3, {\Bbb F}_q)$
such that
\[
{}^t\!A' B_1 = B  {}^t\!A' + \mu_1 B 
           \   \     
         \mbox{\rm  for ${}^t\!A' = \left(
                              \begin{array}{ccc}
                                         0&1&0\\
                                         0&0 &1\\
                                       a&\mu^2 &0
                              \end{array}
                                    \right)$}
\]
when $q=3^e$, because there is a matrix $T \in GL(3, {\Bbb F}_q)$
such that $T^{-1}{}^t\!AT = {}^t\!A' $.
By straightforward computation, we can see that
\[
B = \left(
         \begin{array}{ccc}
               0&1&0\\
               0&\mu  &1\\
               a&2 \mu^2 &2\mu
         \end{array}
     \right)
\]
satisfies the above equation, and $\det B = a \neq 0$.
\qed

\begin{theorem}\label{automgr}
 \begin{enumerate}[{\rm (a)}]
  \item ${\rm Aut}_{{\Bbb F}_q}(C_A)$ contains a cyclic subgroup
   $\langle \overline{B}_0 \rangle$ of order $q^2+q +1$ as a normal subgroup.
  \item Let $ \Lambda_0, \Lambda_1$ and $\Lambda_2$ be
   eigen-vectors of ${}^t\!A$ with distinct eigen-values,
   and $\overline{\Lambda}_0, \overline{\Lambda}_1, \overline{\Lambda}_2$
   the corresponding points of ${\Bbb P}^2$ to these three vectors.
   Then the fixed points of $\overline{B}_0^s$ with $ 1 \leq s < q^2+q+1$
   are $\{ \overline{\Lambda}_0, \overline{\Lambda}_1, 
                                        \overline{\Lambda}_2 \}.$
  \item ${\rm Aut}_{{\Bbb F}_q}(C_A)/\langle \overline{B}_0 \rangle$
  is either trivial or cyclic of order $3$.
  \item ${\rm Aut}_{{\Bbb F}_q}(C_A)/\langle \overline{B}_0 \rangle$
   is nontrivial if and only if either
    \begin{enumerate}[{\rm (i)}]
      \item $q \equiv 1 \mod 3$ and there exists $A' \in  GL(3, {\Bbb F}_q)$
      with $f_{A'}(t) = t^3 - a$ such that $C_A = C_{A'}$, or
      \item $q = 3^e$ and the characteristic polynomial $f_A(t)$
      of $A$ is of the form $f_A(t) = t^3 -(\mu^2 t +a)$
      for some $\mu \in {\Bbb F}_q^{\times}.$
    \end{enumerate}
 \end{enumerate}
\end{theorem}
\proof
(a) This is a special case of Theorem~A.\ref{difference}(b) in Appendix.

(b) Since $B_0^s \in Z_{GL(3, {\Bbb F}_q)}({}^t\!A)$,
three points $\overline{\Lambda}_0, \overline{\Lambda}_1, 
                                        \overline{\Lambda}_2$
are fixed points of $\overline{B}_0^s$ by Lemma~\ref{kerandim}.
Since the eigen-values of $B_0$ are
$\{ \rho, \rho^q, \rho^{q^2} \}$,
where $\overline{\rho}$ is a generator of the cyclic group
${\Bbb F}_{q^3}^{\times}/{\Bbb F}_q^{\times}$
(see the proof of Lemma~A.\ref{centerinGL}),
the eigen-values $\{ \rho^s. \rho^{sq}, \rho^{sq^2} \}$
of $B_0^s$ are distinct each other.
So there is no other fixed point of $\overline{B}_0^s$.

(c) We already saw this in Lemma~\ref{kerandim}.

(d) When we can choose $A' \in GL(3, {\Bbb F}_q)$
so that $C_{A'} = C_A$ and
${\rm Aut}_{{\Bbb F}_q}(C_{A'})
  =Z_{PGL(3, {\Bbb F}_q)}({}^t\!\overline{A'}),$
this is a special case of Theorem~A.\ref{difference} in Appendix.
When we cannot choose such $A'$,
we already saw in Lemma~\ref{replace}.
\qed

\ssgyokan

We can classify each case of (d) in Theorem~\ref{automgr}
up to projective equivalence.

\begin{proposition}\label{classification}
\begin{enumerate}[{\rm (i)}]
 \item When $q \equiv 1 \mod 3$, fix $a \in {\Bbb F}_q^{\times}$
 which is not a cube of any element of ${\Bbb F}_q^{\times}$.
 If ${\rm Aut}_{{\Bbb F}_q}(C_A)/\langle \overline{B}_0 \rangle $
   is nontrivial, then $C_A$ is projectively equivalent to $C_{A'}$
   with
   \[
   A' = \left(
         \begin{array}{ccc}
               0&0&a\\
               1&0 &0\\
               0&1 &0
         \end{array}
     \right)
     \mbox{\ \rm or\ }
       \left(
         \begin{array}{ccc}
               0&0&a^{-1}\\
               1&0 &0\\
               0&1 &0
         \end{array}
     \right).
   \]
 \item When $q = 3^e$, fix $\mu \in {\Bbb F}_q^{\times}$ such that
     $t^3 - (\mu^2 t +1)$ is irreducible.
     If ${\rm Aut}_{{\Bbb F}_q}(C_A)/\langle \overline{B}_0 \rangle $
   is nontrivial, then $C_A$ is projectively equivalent to $C_{A'}$
   with
   \[
   A' = \left(
         \begin{array}{ccc}
               0&0&1\\
               1&0 &\mu^2\\
               0&1 &0
         \end{array}
     \right).
   \]
\end{enumerate}
\end{proposition}
\proof
(i) Let $({\Bbb F}_q^{\times})^3$ be the image of the 3rd power map
${\Bbb F}_q^{\times} \ni \kappa \mapsto \kappa^3 \in {\Bbb F}_q^{\times}.$
Since $q \equiv 1 \mod 3$, the kernel of this map is of order $3$.
Hence $t^3 -a$ $(a \in {\Bbb F}_q^{\times})$
is irreducible over ${\Bbb F}_q$ if and only if
$a \not\in {\Bbb F}_q^{\times} \setminus ({\Bbb F}_q^{\times})^3,$
and ${\Bbb F}_q^{\times}/({\Bbb F}_q^{\times})^3$is of order $3$.
We want to classify the set of polynomials
$\{t^3 -a \mid a \in {\Bbb F}_q^{\times}\setminus ({\Bbb F}_q^{\times})^3 \}$
by the equivalence relation (\ref{equivrelcubics}),
that is,
two monic cubics $f(t)$ and $g(t)$ are equivalent each other
if $f(t) = \rho^3 g(\frac{t-\mu}{\rho})$
for some $\rho \in {\Bbb F}_q^{\times}$ and
$\mu \in {\Bbb F}_q.$
Since ${\Bbb F}_q^{\times}/({\Bbb F}_q^{\times})^3$is of order $3$,
the complete set of representatives of the above set of cubics modulo
the equivalence relation is
$\{ t^3 - a, t^3 - a^{-1} \}.$

(ii) We have to classify the polynomials
$t^3 -(\mu^2 t +a)$ that are irreducible,
by the equivalence relation (\ref{equivrelcubics}).
We show that such polynomials are equivalent one another.
Fix an irreducible polynomial
$f(t) = t^3 - (\mu^2 t + a)$
and choose a root $\lambda$ of $f(t)=0$.
Then $\lambda (\lambda - \mu)(\lambda +\mu) = a$.
Moreover,
two equations
$\lambda +(\lambda - \mu) + (\lambda +\mu) =0$
and
$\lambda (\lambda - \mu) + (\lambda - \mu)(\lambda +\mu)
   + (\lambda +\mu)\lambda = -\mu^2$
are automatically,
because $q = 3^e$.
Hence the three roots of $f(t)=0$ are
$\lambda$, $\lambda - \mu$ and $\lambda +\mu.$
Since the coefficients of $f(t)$ are in ${\Bbb F}_q$,
those three roots coincide with
$\{ \lambda , \lambda^q, \lambda^{q^2} \}.$
So we have $\lambda^q = \lambda + \mu$,
after changing the sign of $\mu$ if need be.
Choose another irreducible polynomial
$g(t) = t^3 - ({\mu'}^2 t + a')$
and a root $\lambda'$ of $g(t)=0.$
Then ${\lambda'}^q = \lambda' + \mu'$ also holds.
Since $1, \lambda, \lambda^2$ form a basis of
${\Bbb F}_{q^3} = {\Bbb F}_q[\lambda] = {\Bbb F}_q[\lambda']$,
there are $\alpha_0, \alpha_1, \alpha_2 \in {\Bbb F}_q$
such that
$\lambda' = \alpha_0  + \alpha_1 \lambda + \alpha_2 \lambda^2$
and either $\alpha_1$ or $\alpha_2$ is a nonzero.
Then we have 
\begin{eqnarray*}
{\lambda'}^q &=& \alpha_0 + \alpha_1 (\lambda + \mu) 
                 + \alpha_2 (\lambda + \mu)^2 \\
             &=& (\alpha_0 + \alpha_1 \mu+ \alpha_2 \mu^2)
                + (\alpha_1 + 2\mu \alpha_2)\lambda + \alpha_2\lambda^2
\end{eqnarray*}
and
\[
{\lambda'}^q = \lambda' + \mu'
   = (\alpha_0 + \mu') + \alpha_1 \lambda + \alpha_2 \lambda^2.
\]
Hence $\alpha_2 =0$, and hence $\lambda' = \alpha_0 + \alpha_1 \lambda.$
So $\lambda'$ is a root of $f( \frac{t - \alpha_0}{\alpha_1} )=0.$
Therefore the monic cubic
$\alpha_1^3 f( \frac{t - \alpha_0}{\alpha_1} )$
must coincide with the minimal polynomial $g(t)$ of $\lambda'$
over ${\Bbb F}_q$.
Moreover, we can find $\rho \in {\Bbb F}_q^{\times}$
so that $\rho^3 = a^{-1}$
because $q=3^e$.
So $\rho^3 f(\frac{t}{\rho}) = t^3 -((\mu\rho)^2t +1).$
This completes the proof of the uniqueness.

Finally we show the existence of an irreducible polynomial
of the form $t^3 -(\mu^2t +a).$
Fix an element $\mu \in {\Bbb F}_q^{\times}$ and consider
the map $\varphi_{\mu}: {\Bbb F}_q \to {\Bbb F}_q$
defined by $\varphi_{\mu}(u) = u(u + \mu)(u-\mu)$.
Then the cardinality of ${\rm Im}\, \varphi_{\mu}$
is $q/3$. Hence ${\Bbb F}_q \setminus {\rm Im}\, \varphi_{\mu}$
is nonempty, and
if we choose $a$ in this set, then the polynomial is irreducible.
\qed

\ssgyokan

\begin{remark}
The cases (i) and (ii) in Proposition~\ref{classification}
correspond to the cases of {\em equianharmonic} and {\em harmonic}
in the sense of Tallini \cite{tal2}, respectively.
He claimed the number of ${\rm Aut}_{{\Bbb F}_q}(C_A)$ was
$6(q^2+q+1)$ if $C_A$ was harmonic (loc. cit., p. 460),
but actually it is $3(q^2+q+1)$,
as was shown in Theorem~\ref{automgr}.
\end{remark}

\ssgyokan

As was pointed out by an anonymous reviewer,
the irreducibility of $C_A$ implies its smoothness
as follows.

\ssgyokan

\noindent
{\em  Another proof of Theorem~{\rm \ref{mainth}},
under the irreducibility of $C_A$.}
Suppose $C_A$ has a singular point $R$ other than
$\overline{\Lambda}_0,
\overline{\Lambda}_1,
\overline{\Lambda}_2 .$
Then from (a) and (b) of Theorem~\ref{automgr},
$C_A$ has at least $q^2+q+1$ singular points
$\{ \overline{B_0}^sR \mid 0 \leq s < q^2+q+1 \}.$
Since the number of singular points of an irreducible curve is
at most the arithmetic genus of the curve,
$q^2 + q+1 \leq \frac{1}{2}(q+1)q$,
which is impossible.
Hence the possibilities of a singular point are only
$\overline{\Lambda}_0$ or
$\overline{\Lambda}_1$ or
$\overline{\Lambda}_2 .$
If we choose a canonical form of the curve as
\[
A= \left(
         \begin{array}{ccc}
               0&0&a\\
               1&0 &b\\
               0&1 &c
         \end{array}
     \right),
\]
these three points are the three points
$\{ Q_{\lambda} \mid
\mbox{\rm $\lambda$ is a root of $f_A(t)=0$}\}$
appeared in the proof of Theorem~\ref{mainth}.
As the computation (\ref{tanCatQ}) in the proof,
$Q_{\lambda}$ is a nonsingular point of $C_A$.
\qed

\gyokan

\noindent
{\bf\Large Appendix}

\sgyokan

\noindent
Throughout Appendix, we fix a matrix $A_0 \in GL(n,{\Bbb F}_q)$
whose characteristic polynomial
$
f_{A_0}(t) = t^n - (a_{n-1}t^{n-1} + \ldots + a_1t + a_0)
$
is irreducible over ${\Bbb F}_q$.
\begin{alemma}\label{centerinGL}
The center $Z_{ GL(n,{\Bbb F}_q)}(A_0)$ of
$A_0 \in GL(n,{\Bbb F}_q)$
is a cyclic group of order $q^n -1$.
\end{alemma}
\proof
Since $f_{A_0}(t)$ is irreducible over ${\Bbb F}_q$,
we may assume that
\begin{equation}\label{generalcanonical}
A_0 = 
  \left(
         \begin{array}{ccccc}
               0&1&0& \cdots & 0 \\
               0&0 &1& \cdots &0 \\
                &  &  & \ddots &\vdots \\
               0&0&0 & \cdots &1\\
              a_0&a_1&a_2&\cdots&a_{n-1}
         \end{array}
     \right).
\end{equation}
Let $\lambda$ be a root of $f_{A_0}(t)=0$.
Then its all roots are
$\{ \lambda, \lambda^q , \ldots , \lambda^{q^{n-1}} \}.$
It is easy to see that
$\Lambda_i = {}^t\! ( 1, \lambda^{q^i}, \lambda^{2q^i},
      \ldots , \lambda^{(n-1)q^i} )$
is an eigen-vector of $A_0$ with the eigen-value $\lambda^{q^i}$
for $i= 0, 1, \ldots , n-1.$
Let
$\Lambda = (\Lambda_0, \Lambda_i, \cdots , \Lambda_{n-1})$.
Then
\[
\Lambda^{-1} A_0 \Lambda
  =  \left(
       \begin{array}{cccc}
         \lambda &         &   & \\
                 &\lambda^q&   & \\
                 &        & \ddots& \\
                 &        &   & \lambda^{q^{n-1}}
       \end{array}
     \right).
\]
Since $\lambda, \lambda^q, \ldots , \lambda^{q^{n-1}} \in {\Bbb F}_{q^{n}}$
are distinct from one another,
\[
Z_{ GL(n,{\Bbb F}_{q^n})}(\Lambda^{-1} A_0 \Lambda)
= \left\{
\left(
       \begin{array}{ccc}
         
                 \beta_1&   & \\
                         & \ddots& \\
                         &   & \beta_n
       \end{array}
     \right)  \mid 
     \beta_1, \ldots , \beta_n \in {\Bbb F}_{q^n}^{\times}
\right\}.
\]
Hence
\[
Z_{ GL(n,{\Bbb F}_q)}(A_0)
= \left\{
\Lambda \left(
       \begin{array}{ccc}
         
                 \beta_1&   & \\
                         & \ddots& \\
                         &   & \beta_n
       \end{array}
     \right) \Lambda^{-1} \mid 
     \beta_1, \ldots , \beta_n \in {\Bbb F}_{q^n}^{\times}
\right\}
\bigcap
GL(n,{\Bbb F}_q).
\]
Let
\[
P=
\left(
         \begin{array}{ccccc}
               0&0& \cdots & 0 &1 \\
               1&0 & \cdots &0 & 0\\
               0&1 & \cdots &0 & 0\\
             \vdots&\vdots &\ddots &\vdots&\vdots \\
               0&0& \cdots&1 &0
         \end{array}
     \right),
\]
which gives the cyclic permutation to the $n$ rows as
$(1, 2, \ldots, n)$ and that to the $n$ columns as
$(n, n-1, \ldots , 1)$.
Hence
$ \Lambda^{(q)} = (\Lambda_1, \Lambda_2, \cdots, \Lambda_{n-1}, \Lambda_0)
   = \Lambda P,$
because $\Lambda_i^{(q)} = \Lambda_{i+1}$ $(i=0, \ldots , n-2)$
and $\Lambda_{n-1}^{(q)} = \Lambda_0$.
We show that
for
\[
B = \Lambda \left(
       \begin{array}{ccc}
         
                 \beta_1&   & \\
                         & \ddots& \\
                         &   & \beta_n
       \end{array}
     \right) \Lambda^{-1} \in Z_{ GL(n,{\Bbb F}_{q^n})}(A_0),
\]
it is in $Z_{ GL(n,{\Bbb F}_q)}(A_0)$
if and only if
\begin{equation}\label{diagonalbeta}
\beta_1^q = \beta_2, \, \beta_2^q = \beta_3,
\ldots, \beta_{n-1}^q = \beta_n, \, \beta_n^q = \beta_1.
\end{equation}
In fact, $B \in GL(n,{\Bbb F}_q)$ if and only if
$B^{(q)} =B$.
Since
\begin{eqnarray*}
 B^{(q)} &=& \Lambda^{(q)} \left(
       \begin{array}{ccc}
         
                 \beta_1^q&   & \\
                         & \ddots& \\
                         &   & \beta_n^q
       \end{array}
     \right) (\Lambda^{-1})^{(q)} \\
     &=&\Lambda P\left(
       \begin{array}{ccc}
         
                 \beta_1^q&   & \\
                         & \ddots& \\
                         &   & \beta_n^q
       \end{array}
     \right) P^{-1}\Lambda^{-1}
     =\Lambda \left(
       \begin{array}{cccc}
         
                 \beta_n^q&   &  &\\
                          &\beta_1^q& & \\
                        &          & \ddots& \\
                         &   & &\beta_{n-1}^q
       \end{array}
     \right) \Lambda^{-1},
\end{eqnarray*}
$B^{(q)} =B$ if and only if the condition (\ref{diagonalbeta})
holds true.

Now we choose a primitive element $\rho$ of ${\Bbb F}_{q^n}.$
Put
\[
B_0 =
\Lambda \left(
       \begin{array}{cccc}
         
                 \rho&   &  &\\
                          &\rho^q& & \\
                        &          & \ddots& \\
                         &   & &\rho^{q^{n-1}}
       \end{array}
     \right) \Lambda^{-1}.
\]
Then $B_0 \in Z_{ GL(n,{\Bbb F}_q)}(A_0)$
because of the condition (\ref{diagonalbeta}).
By definition, $\Lambda_0$ is an eigen-vector of $B_0$
with the eigen-value $\rho$.
Hence the order of $B_0$ is just $q^n-1$.
Moreover, for any $B \in Z_{ GL(n,{\Bbb F}_q)}(A_0)$,
we can find an integer $s$ such that
$B=B_0^s$ because of (\ref{diagonalbeta}) again.
This completes the proof.
\qed

\ssgyokan

Considering the natural homomorphism
$
GL(n,{\Bbb F}_q) \to PGL(n,{\Bbb F}_q),
$
we indicate the image of an object in $GL(n, {\Bbb F}_q)$
by a bar over the object.
\begin{atheorem}\label{difference}
 \begin{enumerate}[{\rm (a)}]
  \item There is a canonical exact sequence{\rm :}
  \[
  0 \to \overline{Z_{ GL(n,{\Bbb F}_q)}(A_0)}
   \to Z_{ PGL(n,{\Bbb F}_q)}(\overline{A}_0)
     \stackrel{\pi}{\to} S_n.
  \]
  \item $\overline{Z_{ GL(n,{\Bbb F}_q)}(A_0)}$ is a cyclic group
  of order $(q^n -1)/(q-1)$.
  \item ${\rm Im}\, \pi$ is a cyclic group,
  which may be trivial.
  \item ${\rm Im}\, \pi$ is nontrivial if and only if
  there is a natural number $k$ with $k >1$ and $k | n$
  such that $q \equiv 1 \mod k$ and
  the characteristic polynomial of $A_0$ is of the form
  $f_{A_0}(t) = 
    t^n - (\sum_{\nu = 0}^{\frac{n}{k}-1}a_{k\nu}t^{k\nu}).$
 \end{enumerate}
\end{atheorem}
\proof
First we define the homomorphism $\pi$.
Let $\overline{B} \in Z_{ PGL(n,{\Bbb F}_q)}(\overline{A}_0)$.
$A_0 B = \rho BA_0$  for some $\rho \in {\Bbb F}_q^{\times}$.
Hence for an eigen-vector $\Lambda_i$ of $A_0$
with the eigen-value $\lambda^{q^i}$,
we have
$A_0 B \Lambda_i = \rho B A_0 \Lambda_i = \rho\lambda^{q^i}B\Lambda_i $,
which means $B \Lambda_i$ is also
an eigen-vector of $A_0$.
Therefore $\overline{B}$ induces a permutation
of $n$ points
$\{ \overline{\Lambda}_0, \overline{\Lambda}_1, \ldots ,
\overline{\Lambda}_{n-1}\} \subset {\Bbb P}^{n-1}({\Bbb F}_{q^n}).$
This is the definition of $\pi$.

The proof of the exactness of the sequence in (a) is similar to
that of Lemma~\ref{kerandim}. So we omit it.

(b) is a consequence of Lemma~A.\ref{centerinGL} with
$Z_{GL(n, {\Bbb F}_q)}(A_0) \supset 
     \{ \kappa E \mid \kappa \in {\Bbb F}_q^{\times} \}.$

Next we prove (c). For two matrices
$B, B' \in GL(n,{\Bbb F}_q)$
with $A_0 B = \rho B A_0$ and $A_0 B' = \rho' B' A_0$,
if $\overline{B} = \overline{B}'$,
then $\rho = \rho'$.
Hence we have another exact sequence
\[
0 \to \overline{Z_{ GL(n,{\Bbb F}_q)}(A_0)}
   \to Z_{ PGL(n,{\Bbb F}_q)}(\overline{A}_0)
     \stackrel{\pi'}{\to} {\Bbb F}_q^{\times}
\]
by $\pi'(\overline{B}) = \rho.$
Since any subgroup of ${\Bbb F}_q^{\times}$ is cyclic,
so is ${\rm Im}\, \pi' \simeq {\rm Im}\, \pi.$

Lastly, we prove (d).
Suppose ${\rm Im}\, \pi$ is nontrivial.
Let $k$ be the order of ${\rm Im}\, \pi \simeq {\rm Im}\, \pi'.$
Then $k>1$ and $k | q-1$.
Let $\pi(\overline{B})$ be a generator of the cyclic group ${\rm Im}\, \pi$.
Note that if $\pi(\overline{B}')$ has a fixed point $\overline{\Lambda}_i$,
then $\pi(\overline{B}')$ is the identity in $S_n$,
which can be shown by a standard argument
similar to the proof of Lemma~\ref{kerandim} (a).
Since $\pi(\overline{B})$ is of order $k$ and
$\pi(\overline{B})^s$ has no fixed point for any $s$ with $1 \leq s <k$,
\[
\pi(\overline{B}) =
(i_{11},\ldots,i_{1k})(i_{21},\ldots,i_{2k})\cdots (i_{l1}, \ldots, i_{lk}),
\]
where the $(i_{i1}, \ldots, i_{ik})$'s
are cyclic permutations of length $k$
and
\[
\{i_{11},\ldots,i_{1k};i_{21},\ldots,i_{2k};\ldots 
    ;i_{l1}, \ldots, i_{lk} \} = \{ 1, 2, \ldots , n\}.
\]
In particular, $k|n$.
Let $\rho = \pi'(\overline{B}) \in {\Bbb F}_q^{\times}.$
Then $\rho$ is a primitive $k$th root of $1$.
Since
$B^{-1}A_0B = \rho A_0$,
$f_{A_0}(t) = f_{\rho A_0}(t) = \rho^n f_{A_0}(\frac{t}{\rho}).$
Hence
\[
f_{A_0}(t) = t^n-(a_{n-1}t^{n-1} + \ldots +a_0)
 = t^n-(a_{n-1}\rho t^{n-1} +a_{n-2}\rho^2 t^{n-2}+ \ldots +a_0\rho^n).
\]
Therefore $a_{n-j} =0$ if $j \not\equiv 0 \mod k$,
that is,
$f_{A_0}(t) = 
    t^n - (\sum_{\nu = 0}^{\frac{n}{k}-1}a_{k\nu}t^{k\nu}).$

Conversely, suppose there is a natural number $k$ with
$k >1$, $k|n$, $q \equiv 1 \mod k$, and
$a_{n-j} =0$ if $j \not\equiv 0 \mod k$.
We may assume that $A_0$ is of the form (\ref{generalcanonical}).
Let $\rho$ be a primitive $k$th root of $1$ and
\[
B = \left(
       \begin{array}{cccc}
         
                 1&   &  &\\
                          &\rho& & \\
                        &          & \ddots& \\
                         &   & &\rho^{n-1}
       \end{array}
     \right).
\]
Then $A_0B = \rho B A_0$.
This completes the proof.
\qed


\end{document}